\documentclass[12pt]{article}
\usepackage{times}
\usepackage{amsmath,amsxtra,amsthm,amssymb,makeidx,graphics,hyperref,graphicx}
\usepackage{latexsym,color}
\theoremstyle{plain}
\newtheorem{theorem}{Theorem}[section]
\newtheorem{lemma}[theorem]{Lemma}
\newtheorem{corollary}[theorem]{Corollary}

\theoremstyle{definition}
\newtheorem{remark}[theorem]{Remark}
\newtheorem{definition}[theorem]{Definition}
\newtheorem{example}[theorem]{Example}

\newtheorem{exercise}[theorem]{Exercise}

\theoremstyle{remark}

\newcommand{\ncom}{\newcommand}

\newcommand{\cA}{{\cal A}}
\newcommand{\cB}{{\cal B}}

\newcommand{\cE}{{\cal E}}

\newcommand{\cI}{{\cal I}}

\newcommand{\cN}{{\cal N}}

\newcommand{\C}{{\mathbb C}}
\newcommand{\Fq}{{{\mathbb F}_q}}
\newcommand{\F}{{{\mathbb F}}}
\newcommand{\bqo}{{{B_q^\partial}}}
\newcommand{\bo}{{{B^\partial}}}

\newcommand{\qb}[2]{{{ {{#1}\choose {#2}}_q }}}
\newcommand{\bin}[2]{{{ {#1}\choose {#2}}} }

\ncom{\noi}{\noindent}
\ncom{\bq}{\begin{equation}}
\ncom{\eq}{\end{equation}}  
\ncom{\beqn}{\begin{eqnarray*}}
\ncom{\eeqn}{\end{eqnarray*}}  
\ncom{\ba}{\begin{array}}
\ncom{\ul}{\underline}   
\ncom{\vphi}{\varphi}   

\ncom{\ea}{\end{array}}
\ncom{\beq}{\begin{eqnarray}}
\ncom{\eeq}{\end{eqnarray}}  
\ncom{\rar}{\rightarrow}
\ncom{\lrar}{\longrightarrow}
\ncom{\bc}{\begin{center}}
\ncom{\ec}{\end{center}}  

\ncom{\bt}{\begin{theorem}}
\ncom{\bcon}{\begin{con}}
\ncom{\et}{\end{theorem}}
\ncom{\econ}{\end{con}}
\ncom{\bl}{\begin{lemma}}
\ncom{\el}{\end{lemma}}  
\ncom{\bco}{\begin{corollary}} 
\ncom{\ds}{\displaystyle}
\ncom{\eco}{\end{corollary}}   
\ncom{\bp}{\begin{pro}}  
\ncom{\ep}{\end{pro}}    
\ncom{\bex}{\begin{example}}
\ncom{\eex}{\end{example}}  
\ncom{\bexr}{\begin{exercise}}
\ncom{\eexr}{\end{exercise}}  
\ncom{\bprob}{\begin{problems}}
\ncom{\eprob}{\end{problems}}  
\ncom{\bd}{\begin{definition}}
\ncom{\ed}{\end{definition}}  
\ncom{\brm}{\begin{remark}}   
\ncom{\erm}{\end{remark}}     
\ncom{\ol}{\overline}
\ncom{\wh}{\widehat} 

\ncom{\pf}{\noi {\bf Proof.  }}
\ncom{\eprf}{\noi {$\Box$}}
\ncom{\be}{\begin{enumerate}} 
\ncom{\ee}{\end{enumerate}}   
\ncom{\seq}{\subseteq}
\newcommand{\inp}[2]{\langle {#1} ,\,{#2} \rangle}

\begin{document}

\title{\bf{{{Eigenbasis for a weighted adjacency matrix associated with the
projective geometry $B_q(n)$}}}}

\author{
{\bf {Murali K. Srinivasan}} \\
{\em  {Department of Mathematics}}\\
{\em  {Indian Institute of Technology, Bombay}}\\
{\em  {Powai, Mumbai 400076, INDIA}}\\
{\bf  \texttt{murali.k.srinivasan@gmail.com}}}
 
\date{}
\maketitle
\begin{center}
{\small{\em To the memory of Reddy}}
\end{center}

\begin{abstract} 

In a recent article {\em Projective geometries, $Q$-polynomial structures,
and quantum groups} Terwilliger (arXiv:2407.14964) defined a certain weighted
adjacency matrix, depending on a free (positive real) parameter, 
associated with the projective geometry, and showed (among many other
results) that it
is diagonalizable, with the eigenvalues and their multiplicities
explicitly written down, and that it satisfies the $Q$-polynomial property
(with respect to the zero subspace).  
 
In this note we
\begin{itemize}
\item Write down an explicit eigenbasis for this matrix. 

\item Evaluate the adjacency matrix-eigenvector products, 
yielding a new proof for the
eigenvalues and their multiplicities. 

\item Evaluate the dual adjacency matrix-eigenvector products 
and directly show that the action of the dual adjacency matrix on the eigenspaces 
of the adjacency matrix is block-tridiagonal, 
yielding a new proof of the $Q$-polynomial property.

\end{itemize}
\end{abstract}
     
\noi

\section{{Introduction}}

Let $q$ be a prime power and let $B_q(n)$ denote the {\em projective
geometry}, i.e., the poset (under inclusion) 
of all subspaces of  $\F_q^n$, 
the $n$-dimensional $\Fq$-vector space of all column vectors 
with $n$ components. 
The set of $k$-dimensional subspaces in $B_q(n)$ is denoted $B_q(n,k)$
and its cardinality is the {\em $q$-binomial coefficient $\qb{n}{k}$}.  
The {\em Galois number}
$$G_q(n)=\sum_{k=0}^n \qb{n}{k}$$
is the total number of subspaces in $B_q(n)$. 

Let $\vphi$ be a positive real number. In {\bf\cite{t5}}, Terwilliger defined
the following (complex) $B_q(n)\times B_q(n)$ {\em{weighted adjacency
matrix}} $A_n$  with entries given by
        \begin{eqnarray}\label{eq:mt1}
        A_n(X,Y) &=& \left\{ \ba{ll}  
        1 & \mbox{if }Y\seq X, \dim(X)= \dim(Y)+1,\\
        \vphi q^{\dim(X)} & \mbox{if }X\seq Y, \dim(Y)=\dim(X)+1, \\
        \frac{\vphi -1}{q-1} q^{\dim(X)} & \mbox{if }X = Y, \\
        0 & \mbox{otherwise.} \ea\right.
        \end{eqnarray}
Terwilliger {\bf\cite{t5}} develops a beautiful theory of these
matrices (generalizing the $\vphi=1$ case studied in {\bf\cite{t3}}. Note
that in this case the diagonal entries are 0). 
The topics studied include: diagonalizability, eigenvalues and
their multiplicities, dual matrix and the (generalized)
$Q$-polynomial property, interpretation of
the $Q$-polynomial structure using the quantum group
$U_{q^{1/2}}(\mathfrak{sl_2})$ in the equitable presentation, and 
split decompositions in the new $Q$-polynomial structure. Here we consider
only the following two basic results.

\bt \label{ev} 
{\bf{(Terwilliger {\bf\cite{t5}})}} $A_n$ is diagonalizable. 
The eigenvalues of $A_n$ are
$$\frac{\vphi q^{n-k} - q^k}{q-1} \mbox{ with multiplicity }
\qb{n}{k},\;\;k=0,1,\ldots ,n.$$
\et 

Let us make a few informal remarks about the proof of Theorem \ref{ev} in
{\bf\cite{t5}}.  It proceeds by considering the irreducible modules
for the Terwilliger algebra $T_n$ generated by $A_n$ and the dual adjacency
matrix $A_n^*$ (defined below). The decomposition of the 
standard module (the complex vector space of column vectors with components
indexed by $B_q(n)$) into  an orthogonal (under the standard
inner product) direct sum of irreducible $T_n$-submodules is known 
({\bf\cite{t1}}). With respect to a suitable basis of an irreducible
$T_n$-submodule $W$,
the matrix representing the restriction of $A_n$ to $W$ is also known
{\bf\cite{t1}} and (a diagonally similar matrix) 
appears in Terwilliger's Leonard
system classification {\bf\cite{t2}} and as such its eigenvalues were known. 
Coupled
with the known decomposition into irreducibles mentioned above 
this yields the proof of
Theorem \ref{ev}.  As pointed out by Terwilliger (in a personal
communication) the eigenvectors for the restriction of $A_n$ to $W$ can also
be written down with components given by certain $_3\phi_2$ hypergeometric
series. However this procedure does not give us an attractive (global)
eigenbasis for $A_n$ itself.

Let $E_n(k),\; k=0,1,\ldots ,n$ denote the eigenspace of $A_n$ corresponding
to the eigenvalue $\frac{\vphi q^{n-k} - q^k}{q-1}$. 
Set $E_n(-1)=E_n(n+1)=\bf{0}$ (the zero subspace).

In {\bf\cite{t5}} Terwilliger defines the {\em{dual adjacency matrix}}
$A_n^*$, a
diagonal $B_q(n)\times B_q(n)$ matrix,  
with diagonal entries given by
$$A_n^*(X,X)=q^{-\dim(X)},\;X\in B_q(n).$$

Note that the eigenvalues of $A_n^*$ are $q^{-k}$ with multiplicities
$\qb{n}{k},\;k=0,1,\ldots ,n$. Also note that $A_n$ acts on the eigenspaces
of $A_n^*$ in a block tridiagonal fashion.

Terwilliger {\bf\cite{t4}} (see Section 20) initiated the extension of the
$Q$-polynomial property to graphs that are not necessarily distance-regular
and gave the first attractive example {\bf\cite{t3}} of this concept.
In the context of the present paper, and the matrix $A_n$,
the (generalized) $Q$-polynomial property with respect to $\bf{0}$
says the following.
\bt \label{qp} 
{\bf{(Terwilliger {\bf\cite{t5}})}} 
\beq \label{trid} 
A_n^*(E_n(k)) & \seq & E_n(k-1) \oplus E_n(k) \oplus E_n(k+1)
,\;\;k=0,1,\ldots ,n.
\eeq
\et 

That is, $A_n^*$ acts on the eigenspaces of $A_n$ in a block tridiagonal
way.
The tridiagonal relations satisfied by $A_n,A_n^*$ are first derived in
{\bf\cite{t5}}  and are then used to prove Theorem \ref{qp}.

The purpose of this note is to write down an explicit eigenbasis for $A_n$.
This yields new proofs of Theorems \ref{ev}, \ref{qp}. We show  
that the method given in {\bf\cite{gs}} to write
down an eigenbasis of $A_n$ in the $\vphi =1$ case easily extends to the
general case with a few small changes to account for the parameter $\vphi$.

In a combinatorial study of vector spaces indexing sets for the objects
studied (bases, eigenvectors etc.) play an important role. The
multiplicities in Theorem \ref{ev} suggest that the eigenvectors may be indexed
by subspaces (just like in the classical case of the $n$-cube the
eigenvectors of the adjacency matrix may be indexed by subsets). However,
there does not seem to be any natural way of indexing the eigenvectors by
subspaces. We now describe the indexing set to be used.

The Goldman-Rota recurrence 
for the Galois numbers ({\bf\cite{gr,kc}}) is the identity 
\beq \label{gri}
G_q(n+1) &=& 2G_q(n) + (q^n - 1)G_q(n-1),\;n\geq 1,\;\;
G_q(0)=1,\;G_q(1)=2.\eeq

More generally, we have the following identity, 
\beq \label{nsw}
\qb{n+1}{k} &=& \qb{n}{k} + \qb{n}{k-1} + (q^n - 1)\qb{n-1}{k-1},\;n,k\geq
1,\eeq
with $\qb{0}{k}=\delta(0,k)$ and $\qb{n}{0}=1$.
Note that (\ref{gri}) follows by summing (\ref{nsw}) over $k$.

Let $\cI_q(n)$  
denote the set of all distinct irreducible characters (all of
degree $1$) of the finite abelian group $\F_q^n$
and let $\cN_q(n)$  
denote the set of all distinct nontrivial irreducible characters
of $\F_q^n$. So $|\cN_q(n)|=q^n-1$.

For $n\geq 0$, inductively define a set $\cE_q(n)$ consisting of sequences
as follows (here $()$ denotes the empty sequence):
\beqn
\cE_q(0) &=& \{()\},\\
\cE_q(1) &=& \{(0),(1)\},\\
\cE_q(n) &=& 
\{(\alpha_1,\ldots ,\alpha_t)\;|\;
             (\alpha_1,\ldots \alpha_{t-1})\in
\cE_q(n-1),\;\alpha_t \in \{0,1\}\}\\
&&\cup \{(\alpha_1,\ldots ,\alpha_t)\;|\;
             (\alpha_1,\ldots \alpha_{t-1})\in\cE_q(n-2),\;
\alpha_t \in \cN_q(n-1)\},\;n\geq 2.
\eeqn
Given $\alpha\in \cE_q(n)$, let $N(\alpha)$ denote the number of nonzero
entries in the sequence $\alpha$ ( note that a nonzero entry is either 1 or
an element of $\cN_q(m)$ for some $m$). Set
\beqn
\cE_q(n,k) &=& \{\alpha\in \cE_q(n)\;|\;N(\alpha)=k\},\\
e_q(n,k) &=& |\cE_q(n,k)|.
\eeqn
It is easy to see that
\beq 
e_q(n+1,k) &=& e_q(n,k) + e_q(n,k-1) + (q^n - 1)e_q(n-1,k-1),\;n,k\geq
1,\eeq
with $e_q(0,k)=\delta(0,k)$ and $e_q(n,0)=1$,
the same recurrence (with the same initial conditions) as (\ref{nsw}), Thus
$e_q(n,k)=\qb{n}{k}$ and $|\cE_q(n)|=|B_q(n)|$.

Another relevant structure in the present situation is an inner product
(different from the standard inner product) 
on the (complex) vector space of column vectors
with components indexed by $B_q(n)$. 

For $X\in B_q(n)$ with $\dim(X)=k$ define
$$\pi(X) = \frac{\vphi^{k} q^{\bin{k}{2}}}{P_q(n)},
$$
where $P_q(n) =\prod_{k=0}^{n-1} (1+\vphi q^k)$.
We have
$$\sum_{X\in B_q(n)}\pi(X)\;=\;\frac{\sum_{k=0}^n \vphi^{k} 
q^{\bin{k}{2}}\qb{n}{k}}
{P_q(n)}\;=\;1$$
where the second equality follows by
the $q$-binomial theorem (so  $\pi$ is a probability vector on $B_q(n)$). 
In the statistical literature the distribution $(p_0, p_1,\ldots , p_n)$,
where $p_k = \frac{\vphi^{k} 
q^{\bin{k}{2}}\qb{n}{k}}{P_q(n)}$, is called Kemp's distribution, a two parameter
analog of the binomial distribution (which corresponds to
$\vphi=1$, $q=1$) (see, for example, {\bf\cite{kk}}).

Given vectors $u,v$ with components indexed by $B_q(n)$ define
\beq \label{inp}
\inp{u}{v}_\pi =\sum_{X\in B_q(n)}\ol{u(X)}v(X)\pi(X).
\eeq
Since $P_q(n)$ is independent of $k$, for orthogonality of vectors we do not
need the denominator in the definition of $\pi(X)$. It is there only to make
$\pi$ a probability vector.

It follows from Lemma 10.1 in {\bf\cite{t5}} (also see Section 2 of the
present paper) that $A_n$ is self-adjoint with
respect to the inner product (\ref{inp}).

We now state our two main results. 
\bt \label{mt}
For each $\alpha\in \cE_q(n)$ we define a column 
vector $v_\alpha$, with components indexed by $B_q(n)$,
such that

(i) $A_n v_\alpha = \frac{\vphi q^{n-k} - q^k}{q-1}\;v_\alpha$, 
where $k=N(\alpha)$.

(ii) $\{v_\alpha\;|\; \alpha\in \cE_q(n)\}$ is an eigenbasis of $A_n$.

(iii) $\{v_\alpha\;|\; \alpha\in \cE_q(n)\}$ is orthogonal with respect to
the inner product (\ref{inp}).
\et
Theorem \ref{ev} is an immediate consequence of parts (i) and (ii) of
Theorem \ref{mt}. 

Our second main result, in Section 4, considers 
the action of $A_n^*$ on the
eigenspaces of $A_n$. We evaluate $A_n^* v_\alpha$
and using the structure of the eigenspaces 
developed during the course of proof of Theorem \ref{mt} we directly show
(\ref{trid}), yielding an alternative proof of Theorem \ref{qp}.

Our approach to proving Theorem \ref{mt} and Theorem \ref{qp} is inductive
(which explains our notation $A_n, A_n^*$,  instead of $A,A^*$, as in
{\bf\cite{t5}}). It is based on a decomposition of the space of column
vectors with components indexed by $B_q(n)$ that was worked out in
{\bf\cite{s}} (the indexing set for the eigenbasis is suggested by this
decomposition).

\section{{A decomposition of $\C[B_q(n)]$}}

In this section we recall the decomposition from {\bf\cite{s}}. All the
omitted proofs may be found in Section 2 of {\bf\cite{s}}. This section is
virtually the same as Section 3 of {\bf\cite{gs}} and is included here
primarily for the readers convenience.

Before recalling this decomposition we first observe 
that $A_n$ is self-adjoint under the inner product
(\ref{inp}).
Let $D_n$ be the $B_q(n)\times B_q(n)$
diagonal matrix with diagonal entry in row $X$, column $X$ 
given by $\sqrt{\vphi^kq^{\bin{k}{2}}}$, where $k=\dim(X)$. 
\bl ({\bf\cite{t5}}, Lemma 10.1)
$D_nA_nD_n^{-1}$ is symmetric. Thus $A_n$ is diagonalizable with real
eigenvalues.
\el
\pf
For $X\in B_q(n,k)$,
$Y\in B_q(n,r)$ the entry in row $X$, column $Y$ of $D_nA_nD_n^{-1}$
is given by
        \beqn
        \lefteqn{\sqrt{\vphi^k q^{\bin{k}{2}}}\, A_n(X,Y)\, \sqrt{\vphi^{-r}
                      q^{-\bin{r}{2}}}}\\
        &=&\left\{\ba{ll}
        \sqrt{\phi^{k}q^{\bin{k}{2}}}\, q^k\, 
         \sqrt{\phi^{-(k+1)q^{-\bin{k+1}{2}}}} =
        \sqrt{\vphi q^k} & \mbox{if $X\seq Y$ and $r=k+1$,}\\   
        \sqrt{\vphi^{k}q^{\bin{k}{2}}} \sqrt{\vphi^{-(k-1)}
        q^{-\bin{k-1}{2}}} =  
        \sqrt{\vphi q^{r}} & \mbox{if $Y\seq X$ and $r=k-1$,}\\
        \frac{\vphi - 1}{q-1}q^k  & \mbox{if $X=Y$,}\\ 
        0 & \mbox{otherwise,} 
        \ea
        \right.\\
        &=&\left\{\ba{ll}
        \sqrt{\vphi q^{\min\{\dim(X),\dim(Y)\}}}&
\mbox{if $X\seq Y$ or $Y\seq X$, and $|\dim(X)-\dim(Y)|=1$,} \\
        \frac{\vphi - 1}{q-1}q^{\dim(X)}  & \mbox{if $X=Y$,}\\ 
        0 & \mbox{otherwise,} 
        \ea
        \right.
        \eeqn  
yielding a symmetric matrix. \eprf

Since $P_q(n)$ is independent of $k$,
the argument showing that $D_nA_nD_n^{-1}$
is symmetric shows that $A_n$ is
self-adjoint with respect
to the inner product (\ref{inp}).

For a finite set $S$, we denote the complex vector space with $S$
as basis by $\C[S]$. We can think of $\C[S]$ as the space of all (complex)
column vectors with components indexed by $S$. Thus, for $s\in S$, the
vector $s\in \C[S]$ corresponds to the column vector with 1 in component $s$
and 0's elsewhere.

We have (vector space direct sum)
$$\C[B_q(n)]=\C[B_q(n,0)]\oplus 
\C[B_q(n,1)]\oplus \cdots \oplus\C[B_q(n,n)].$$

An element $v\in \C[B_q(n)]$ can be uniquely written as
$v=v_0 + v_1 +\cdots + v_n$, where $v_k\in \C[B_q(n,k)]$ for all $k$.
The $v_k$ are the {\em homogeneous components} of $v$. We say
$v$ is {\em homogeneous} if 
$v\in \C[B_q(n,k)]$ for some $k$. A subspace of $\C[B_q(n)]$ is {\em
homogeneous} if it contains the homogeneous components of all its elements.

Let $0\leq k \leq n$.
Define a linear operator $I_{n,k}:\C[B_q(n)]\rar \C[B_q(n)]$ 
by $I_{n,k}(X)=0$ if $\dim(X)\not=k$ and
$I_{n,k}(X)= X$ if $\dim(X)=k$, $X\in B_q(n)$.
The {\em $k^{th}$ up operator}  $U_{n,k}:\C[B_q(n)]\rar \C[B_q(n)]$ 
is defined, for $X\in B_q(n)$, by $U_{n,k}(X)=0$ if $\dim(X)\not=k$ and
$U_{n,k}(X)= \sum_{Y} Y$,
where the sum is over all $Y\in B_q(n)$ covering $X$, if $\dim(X)=k$.
Similarly we define the {\em $k^{th}$ down operator} 
$D_{n,k}:\C[B_q(n)]\rar \C[B_q(n)]$ (we have $U_{n,n}=D_{n,0}=0$).
Set $U_n=U_{n,0}+U_{n,1}+\cdots +U_{n,n}$ 
and $D_n=D_{n,0}+D_{n,1}+\cdots +D_{n,n}$, called, respectively,
the {\em up and down} operators on $\C[B_q(n)]$.

It is easily seen that $A_n$ is the matrix of the operator
$$\cA_n = U_n +\sum_{k=0}^n \vphi q^{k-1}D_{n,k} + \sum_{k=0}^n \frac{\vphi -
1}{q-1} q^k I_{n,k}$$
with respect to the {\em standard basis} $B_q(n)$. 

Define $\cA_n^* : \C[B_q(n)]\rar \C[B_q(n)]$ by
\beqn   \cA_n^*(X)&=&q^{-\dim(X)} X,\;\;X\in B_q(n).
\eeqn

In Section 4 we shall study the action of $\cA_n^*$ on the eigenspaces of
$\cA_n$. 

The decomposition of $\C[B_q(n)]$ to be studied gives
a linear algebraic interpretation to (\ref{gri}) and 
(\ref{nsw}). Denote
the standard basis vectors of $\F^n_q$ by the column vectors $e_1,\ldots ,e_n$.
We identify $\F_q^k$, for $k<n$, with the subspace of $\F_q^n$ consisting  of
all vectors with the last $n-k$ components zero. 
So $B_q(k)$ consists of all
subspaces of $\F_q^n$ contained in the subspace spanned by $e_1,\ldots ,e_k$.

Define $\bqo(n+1)$ to be the
collection of all subspaces in $B_q(n+1)$ not contained in the hyperplane
$\F_q^{n}$, i.e.,
$$\bqo(n+1) = B_q(n+1) - B_q(n) = \{ X\in B_q(n+1) : X\not\subseteq \F_q^{n} 
\},\;n\geq 0.$$
For $1\leq k \leq n+1$, let $\bqo(n+1,k)$ denote the set of all subspaces in
$\bqo(n+1)$ with dimension $k$. We consider $\bqo(n+1)$ 
as an induced subposet of $B_q(n+1)$.

We have a direct sum decomposition
\beq \label{bod}
\C[B_q(n+1)] = \C[B_q(n)] \oplus \C[\bqo(n+1)].
\eeq
We shall now give a further decomposition of $\C[\bqo(n+1)]$.

Let $H(n+1,\F_q)$ denote the subgroup of $GL(n+1,\F_q)$ consisting of all
matrices of the form
$$ \left[ \ba{cc}
           I & \ba{c} a_1 \\  \cdot \\ \cdot \\  a_n \ea \\
           0 \cdots 0 & 1
          \ea 
   \right],   
$$
where $I$ is the $n\times n$ identity matrix.

The additive abelian group $\F_q^n$ is isomorphic to $H(n+1,\F_q)$ via
$\phi : \F_q^n \rar H(n+1,\F_q)$ given by
$$\phi \left(\left[ \ba{c} a_1 \\ \cdot \\ \cdot \\a_n \ea \right]\right)
\rar
\left[ \ba{cc}
           I & \ba{c} a_1 \\  \cdot \\ \cdot \\  a_n \ea \\
           0 \cdots 0 & 1
          \ea 
   \right].   
$$   

Recall that $\cI_q(n)$  
denotes the set of all distinct irreducible characters 
of the finite abelian group $H(n+1,\F_q)$
and $\cN_q(n)$  
denotes the set of all distinct nontrivial irreducible characters
of $H(n+1,\F_q)$.

There is a natural (left) permutation action
of $H(n+1,\F_q)$ on $\bqo(n+1)$ and $\bqo(n+1,k)$. This permutation action
induces representations of $H(n+1,\F_q)$ on $\C[\bqo(n+1)]$ and
$\C[\bqo(n+1,k)]$.

For $\chi\in \cI_q(n)$,
let $W(\chi)$ (respectively, $W(\chi,k)$) denote the isotypical component of
$\C[\bqo(n+1)]$ (respectively, $\C[\bqo(n+1,k)]$) corresponding to the 
irreducible representation of $H(n+1,\F_q)$ with character $\chi$.
When $\chi$ is the trivial character we denote $W(\chi)$ (respectively,
$W(\chi,k)$) by  $W(0)$ (respectively, $W(0,k)$).
We have the following decompositions, (note that $W(\chi,n+1)$, 
for $\chi\in \cN_q(n)$, is the zero module). 
\beq
W(0) &=& W(0,1) \oplus \cdots \oplus W(0,n+1), \\
W(\chi) &=& W(\chi,1) \oplus \cdots \oplus W(\chi,n),
\;\;\;\chi \in \cN_q(n), \\ \label{bod1}
\C[\bqo(n+1)] &=& W(0) \oplus \left(\oplus_{\chi \in \cN_q(n)} W(\chi)\right).
\eeq

Note that $\C[B_q(n)], \C[\bqo(n+1)], W(0), W(\chi)$ are homogeneous subspaces
of $\C[B_q(n+1)]$.

Now $GL(n+1,\F_q)$ acts on $B_q(n+1)$ and $\C[B_q(n+1)]$ and the action of
$U_{n+1}$ commutes with the action of $GL(n+1,\F_q)$ (and 
hence with the action of $H(n+1,\F_q)$). 
Also, $\C[\bqo(n+1)]$ is clearly closed under $U_{n+1}$. Thus 
\beq \label{d0}
W(0),\;\;W(\chi),\chi\in \cN_q(n) \mbox{ are  $U_{n+1}$-closed}.
\eeq

Define an equivalence 
relation $\sim$ on $\bqo(n+1)$ by $X\sim Y$ iff 
$X\cap \,\F_q^{n}= Y\cap \,\F_q^{n}$.
Denote the
equivalence class of $X\in \bqo(n+1)$ by $[X]$. 
For a subspace $X\in B_q(n)$, define $\wh{X}$ to be the subspace in $\bqo(n+1)$
spanned by $X$ and $e_{n+1}$.

\bl \label{el}
Let $X,Y \in \bqo(n+1)$ and $Z,T\in B_q(n)$. Then

\noi (i) $\mbox{dim}\,(X\cap \,\F_q^{n})=\mbox{dim}\,X - 1$ and 
$\wh{X\cap \,\F_q^{n}}\in [X]$. 

\noi (ii) $Z\leq T$ iff $\wh{Z} \leq \wh{T}$.

\noi (iii) $Y$ covers $X$ iff 

(a) $Y\cap \,\F_q^{n}$ covers $X\cap \,\F_q^{n}$ and

(b) $Y = \mbox{ span}\,((Y\cap \F_q^{n})\cup \,\{v\})$ for any
$v\in X - \F_q^{n}$. 

\noi (iv) The number of subspaces $Z'\in \bqo(n+1)$ with
$Z'\cap \,\F_q^{n} = Z$ is $q^l$, where $l= n -\mbox{dim}\,Z$.
Thus, $|[X]|=q^{n+1-k}$, where $k=\mbox{ dim}\,X$.\;\;\mbox{\eprf}
\el

For $X\in \bqo(n+1)$, let
$G_X \subseteq H(n+1,\F_q)$ denote the stabilizer of $X$.

\bl \label{orl}
Let $X,Y\in \bqo(n+1)$. Then

\noi (i) The orbit of $X$ under the action of $H(n+1,\F_q)$ is $[X]$.

\noi (ii) Suppose $Y$ covers $X$. Then the bipartite graph of the covering
relations between $[Y]$ and $[X]$ is regular with degrees $q$ (on the
$[Y]$ side) and $1$ (on the $[X]$ side). 

\noi (iii) Suppose $X\subseteq Y$. Then $G_X \subseteq G_Y$.\;\;\eprf
\el

Consider $\C[B_q(n+1)]$. For $X\in B_q(n)$ define
\beqn\theta_n(X)= \sum_Y Y,
\eeqn
where the sum is over all $Y\in \bqo(n+1)$ covering $X$. Equivalently, the sum
is over all $Y\in \bqo(n+1)$ with $Y\cap\,\F_q^n = X$, i.e., $Y\in [\wh{X}]$. It
follows from Lemma \ref{orl}(i) that 
\beq \label{d1}
\theta_n : \C[B_q(n)] \rar W(0)
\eeq
is a linear isomorphism.

Combining (\ref{bod}) and (\ref{bod1}) we have the decomposition
\beq \label{d2}
\C[B_q(n+1)] &=& (\C[B_q(n)]\oplus W(0)) \oplus 
\left(\oplus_{\chi \in \cN_q(n)} W(\chi)\right),
\eeq
where, by (\ref{d0}) and (\ref{d1}), 
\beq \label{d3}
\C[B_q(n)]\oplus W(0)\mbox{ is  $U_{n+1}$-closed}.
\eeq

Let $\psi_k$ (respectively, $\psi$) denote the character of the
permutation
representation of $H(n+1,\F_q)$ on $\C[\bqo(n+1,k)]$ (respectively,
$\C[\bqo(n+1)]$) corresponding to the left action. Clearly 
$\psi = \sum_{k=1}^{n+1} \psi_k$.
Below $[,]$
denotes character inner product and the $q$-binomial coefficient $\qb{n}{k}$ is taken to
be zero when $n$ or $k$ is $< 0$.

\bt \label{grl}
\noi (i) Let $\chi \in \cI_q(n)$ be the trivial character. Then
$[\chi , \psi_k ]= \qb{n}{k-1},\;1\leq k \leq n+1.$

\noi (ii) Let $\chi \in \cN_q(n)$. Then
$[\chi , \psi_k ]= \qb{n-1}{k-1},\;1\leq k \leq n+1.$  \eprf
\et
\brm \label{remark1}
Using Theorem \ref{grl}(ii) we see that
\beq \label{grd}
\dim(W(\chi))&=&\sum_{k=1}^{n+1} \qb{n-1}{k-1}=\;G_q(n-1),\;\;\chi\in\cN_q(n).
\eeq 
Now, by taking dimensions on both sides of (\ref{d2}) 
and using (\ref{d1}), (\ref{grd}) we get the Goldman-Rota identity (\ref{gri}).
More generally, by restricting to dimension $k$ on both sides of
(\ref{d2}), we get the identity (\ref{nsw}).
\erm

For $\chi \in \cI_q(n)$, define the following element of the group algebra of
$H(n+1,\F_q)$:
$$p(\chi)=\sum_g \ol{\chi(g)}\,g,$$
where the sum is over all $g\in H(n+1,\F_q)$. For $1\leq k \leq n+1$, the map
\beq \label{proj}
&p(\chi) : \C[\bqo(n+1,k)] \rar \C[\bqo(n+1,k)],&
\eeq
given by $v\mapsto \sum_{g\in H(n+1,\F_q)} \ol{\chi(g)}\,gv$,
is a nonzero multiple of the $H(n+1,\F_q)$-linear projection onto
$W(\chi,k)$. Similarly for $p(\chi) : \C[\bqo(n+1)] \rar \C[\bqo(n+1)]$.

For future reference we record the following
observation:
\beq \label{ds}
p(\chi)(\wh{Y}) \mbox{ and } p(\chi)(\wh{Z}) 
\mbox{ have disjoint supports, for }
Y\not= Z \in B_q(n).
\eeq 

\bl \label{cl} Let $X\in \bqo(n+1)$ and $\chi\in \cI_q(n)$. 
Then $p(\chi)(X) =
0$ iff $\chi : G_X \rar \C^*$ is a nontrivial character of $G_X$. \eprf
\el

\bt \label{f}

(i) Let $\chi\in \cI_q(n),\;X,Y\in \bqo(n+1)$ with $X=hY$ for
some $h\in H(n+1,\F_q)$. Then
$$p(\chi)(X) = \ol{\chi(h^{-1})}\,p(\chi)(Y).$$

\noi (ii) Let $\chi\in \cI_q(n)$. Then
$\{ p(\chi)(\wh{X}) : X\in B_q(n,k-1) \mbox{ with } p(\chi)(\wh{X})
\not= 0\}$ is a basis of $W(\chi,k)$, $1\leq k \leq n+1$.

\noi (iii) Let $\chi\in \cI_q(n)$ and let $X,Y\in B_q(n)$ with $X$ covering $Y$.
$$p(\chi)(\wh{X})\not= 0 \mbox{ implies } p(\chi)(\wh{Y})\not= 0.
\;\;\mbox{\eprf}$$

\et

Let $\chi \in \cN_q(n)$. 
By Theorem \ref{grl}(ii) we have $\dim(W(\chi,n))=1$.
It thus follows by Theorem \ref{f}(ii) and (\ref{ds}) above 
that there is a unique element $X(\chi)\in B_q(n,n-1)$
such that $p(\chi)(\wh{X(\chi)})\not=0$. Moreover, 

\bl \label{cu}
Let $Y\in B_q(n,n-1)$. Then
$$|\{\chi\in\cN_q(n)\;|\; X(\chi)=Y\}| = q-1.\;\;\mbox{ \eprf}$$
\el

\section{{Eigenvectors of $\cA_n$}}

Consider the decomposition
\beq \label{od} 
\C[B_q(n+1)] &=& (\C[B_q(n)]\oplus W(0)) \oplus 
\left(\oplus_{\chi \in \cN_q(n)} W(\chi)\right).
\eeq
We claim that
\beq 
\C[B_q(n)]\oplus W(0)\mbox{ and }W(\chi),\chi\in \cN_q(n) 
 \mbox{ are  $I_{n+1,k}, D_{n+1,k}$-closed}.
\eeq
This can be seen as follows. Since $\C[B_q(n)], W(0), W(\chi)$ are
homogeneous subspaces of $\C[B_q(n+1)]$ it is easily seen that they are
closed under $I_{n+1,k}$. 
Consider the {\em standard inner product} on
$\C[B_q(n+1)]$ (i.e., declare $B_q(n+1)$ to be an orthonormal basis), which
is $GL(n+1, \F_q)$-invariant (and hence $H(n+1, \F_q)$-invariant). 
It follows that
$W(0)$, $W(\chi),\;\chi\in\cN_q(n)$ are orthogonal and hence the decomposition
(\ref{od}) is orthogonal.
Since $D_{n+1}$ is the adjoint of $U_{n+1}$, it now follows from
(\ref{d0}) and (\ref{d3}) that $\C[B_q(n)] \oplus W(0)$ and $W(\chi)$ are
closed under $D_{n+1}$. The claim now follows from homogeneity.

Thus $\C[B_q(n)]\oplus W(0)$ and $W(\chi),\chi\in\cN_q(n)$ are closed
under $\cA_n$.

Let $X\in B_q(n,k)$. By Lemma \ref{el}(iv) we see that $\theta_n(X)$
is a sum of $q^{n-k}$ subspaces in $\bqo(n+1,k+1)$. 
So we can express 
$D_{n+1,k+1}(\theta_n(X))$ in the following form 
\beq \label{dtheta}
D_{n+1,k+1}(\theta_n(X)) = q^{n-k}X + v,
\eeq
for some $v\in \C[\bqo(n+1,k)]$. A little reflection shows that $v\in W(0)$.
Defining 
$$D^{'}_{n+1,k+1}(\theta_n(X))=v$$ 
gives a linear map
\beqn D^{'}_{n+1,k+1}: W(0)\rar W(0),\;\;0\leq k \leq n,
\eeqn
that takes the vector $\theta_n(X)$ for $X\in B_q(n,k)$ to the vector
$v$ above.

Define $\cA_n^{'}: W(0)\rar W(0)$ by
\beqn
\cA_n^{'} = U_{n+1} + \sum_{k=0}^n \vphi q^k D^{'}_{n+1,k+1} 
                       + \sum_{k=0}^n \frac{\vphi - 1}{q-1}\;q^{k+1}
I_{n+1,k+1}.
\eeqn

We have the following relations (the first of which follows from
(\ref{dtheta})):
\beq \label{ind1}
\cA_{n+1}(\theta_n(v))&=& \vphi q^n v 
       + \cA^{'}_n(\theta_n(v)),\;\;v\in \C[B_q(n)],\\ \label{ind2}
\cA_{n+1}(v)&=& \cA_{n}(v) + \theta_n(v),\;\;v\in \C[B_q(n)].
\eeq

We now write down the matrix of $\cA^{'}_n$ with respect to the basis
$\{\theta_n(X)\;|\;X\in B_q(n)\}$ of $W(0)$.

It follows from Lemma \ref{el}(iii) and Lemma \ref{orl}(ii) that
\beqn U_{n+1}(\theta_n(X)) &=& \sum_Y q\,\theta_n(Y)\;=\;q\,\theta_n(U_n(X))
,\;\;X\in B_q(n),
\eeqn
where the middle sum is over all $Y\in B_q(n)$ covering $X$. 
Similarly, it follows that
\beqn q^k D^{'}_{n+1,k+1}(\theta_n(Y))&=& 
\sum_X q \left\{ q^{k-1}\theta_n(X)\right\}\;=\;q\,\theta_n(q^{k-1} D_{n,k}(Y))
,\;\;Y\in B_q(n,k),
\eeqn
where the middle sum is over all $X\in B_q(n)$ covered by $Y$. Clearly,
\beqn q^{k+1} I_{n+1,k+1}(\theta_n(X)) &=& q\,\theta_n(q^k I_{n,k}(X))
,\;\;X\in B_q(n,k).
\eeqn

Thus we see that 
\beq \label{ind}
&\mbox{Matrix of $\cA^{'}_n$ with respect to the basis
$\{\theta_n(X)\;|\;X\in B_q(n)\}$ is $qA_n$.}
\eeq

For a finite vector space $X$ over $\Fq$ we denote by
$B_q(X)$ the set of all subspaces of $X$ 
and we denote by $U_X$ (respectively, $D_X$) the up
operator (respectively, down operator) on $\C[B_q(X)]$.

Let $(V_1,f_1)$ be a pair consisting of a finite dimensional vector space
$V_1$ (over $\C$) and a linear operator $f_1$ on $V$. Let $(V_2,f_2)$ be
another such pair. By an isomorphism of pairs $(V_1,f_1)$ and $(V_2,f_2)$ 
we mean a linear isomorphism $\tau : V_1
\rar V_2$ such that $\tau(f_1(v)) = f_2(\tau(v)),\;v\in V_1$.

\bt \label{grind} 
Let $\chi\in\cN_q(n)$ and $X=X(\chi)$. Define
\beqn
\lambda(\chi): \C[B_q(X)]\rar W(\chi)
\eeqn
by $Y\mapsto q^{-\dim(Y)}p(\chi)(\wh{Y}),\;\;Y\in B_q(X)$.

Then 

(i) $\lambda(\chi)$ is an isomorphism of pairs $(\C[B_q(X)], qU_X)$
and $(W(\chi),U_{n+1})$.

(ii) $\lambda(\chi)$ is an isomorphism of pairs $(\C[B_q(X)], D_X)$
and $(W(\chi),D_{n+1})$.
\et
\pf This is Theorem 4.2 in {\bf\cite{gs}}. \eprf

Before proceeding further we introduce some notation.
Let
$X\in B_q(n,n-1)$. The pairs $(\C[B_q(X)],U_X)$ and $(\C[B_q(n-1)],U_{n-1})$ 
are clearly isomorphic with many possible isomorphisms. We now define a
canonical isomorphism, based on the concept of a matrix in Schubert normal
form.

A $n\times k$ matrix $M$ over $\F_q$ is in {\em Schubert normal form}
 (or, {\em column reduced echelon form}) provided

\noi (i) Every column is nonzero.

\noi (ii) The last  nonzero entry in every column is a $1$. Let the last
nonzero entry in column $j$ occur in row $r_j$.

\noi (iii) We have $r_1 < r_2 < \cdots < r_k$ and the submatrix of $M$
formed by the rows $r_1,r_2,\ldots ,r_k$ is the $k\times k$ identity matrix.

It is well known
that every $k$ dimensional subspace
of $\F^n_q$ is the column space of a unique $n\times k$ matrix in Schubert
normal form.

Let $X\in B_q(n,n-1)$ and let $M(X)$ be the $n\times (n-1)$ matrix in
Schubert normal form with column space $X$. The map $\tau(X) : \F^{n-1}_q
\rar X$ given by $e_j \mapsto \mbox{ column $j$ of }M(X)$ is clearly a linear
isomorphism and this isomorphism gives rise to an isomorphism
$$\mu(X) : \C[B_q(n-1)]\rar \C[B_q(X)]$$
of pairs $(\C[B_q(n-1)],U_{n-1})$ and $(\C[B_q(X)],U_X)$ 
(and also of pairs $(\C[B_q(n-1)],D_{n-1})$ and $(\C[B_q(X)],D_X)$) 
given by
$\mu(X)(Y)=\tau(X)(Y),\;Y\in B_q(n-1)$. 

Since $\dim(\wh{Y})=1+\dim(Y)$, the following result 
now follows from Theorem \ref{grind}.

\bt \label{grind1}
Let $\chi\in \cN_q(n)$ and $X=X(\chi)$. Then the composition
$\lambda(\chi)\mu(X)$ is an isomorphism of pairs $(\C[B_q(n-1)],q\cA_{n-1})$
and $(W(\chi),\cA_{n+1})$. \eprf
\et

\bl \label{kc}
Let $v\in \C[B_q(n)]$ satisfy
$$
\cA_n(v)= \frac{\vphi q^{n-k} - q^k}{q-1}\;v.$$
Then
\beqn 
\cA_{n+1}(q^k v + \theta_n(v)) &=& \frac{\vphi q^{n+1-k} - q^k}{q-1}\;
         (q^k v + \theta_n(v)),\\
\cA_{n+1}(\vphi q^{n-k} v - \theta_n(v)) &=& \frac{\vphi q^{n+1-(k+1)} -
q^{k+1}}{q-1}\;(\vphi q^{n-k} v - \theta_n(v)).\\
\eeqn
\el
\pf
By (\ref{ind}) we have 
$$
\cA^{'}_n(\theta_n(v))=  \frac{q(\vphi q^{n-k} - q^k)}{q-1}\;v.
$$ 
We have, by (\ref{ind1}) and (\ref{ind2}),
\beqn
\cA_{n+1}(q^k v + \theta_n(v))
                    &=& q^k \cA_{n+1}(v) 
                    + \cA_{n+1}(\theta_n(v))\\
                    &=& q^k(\cA_n(v) +\theta_n(v))
                        +\vphi q^n v + \cA^{'}_n(\theta_n(v))\\
                    &=& \left(\vphi q^{n-k} + \frac{\vphi q^{n-k} -
                    q^k}{q-1}\right)
                    q^k v
                     + \left(q^k +\frac{q(\vphi q^{n-k} - q^k)}{q-1} \right)
                    \theta_n(v)\\
                    &=& \frac{\vphi q^{n+1-k} - q^k}{q-1}(q^k v + 
                    \theta_n(v)),
\eeqn
and
\beqn
\cA_{n+1}(\vphi q^{n-k} v - \theta_n(v))
                    &=& \vphi q^{n-k}\cA_{n+1}(v) 
                    - \cA_{n+1}(\theta_n(v))\\
                    &=& \vphi q^{n-k}(\cA_n(v) +\theta_n(v))
                        -  \vphi q^n v - \cA^{'}_n(\theta_n(v))\\
                    &=& \left(-q^{k} + \frac{\vphi q^{n-k} -
                    q^k}{q-1}\right)
                    \vphi q^{n-k} v
                     - \left(-\vphi q^{n-k} +\frac{q(\vphi q^{n-k} - q^k)}{q-1} \right)
                    \theta_n(v)\\
                    &=& \frac{\vphi q^{n+1-(k+1)} - q^{k+1}}{q-1}
                    (\vphi q^{n-k} v - \theta_n(v)).\;\;\;\;\;\mbox{\eprf}
\eeqn

\noi {\bf{Proof of Theorem \ref{mt} (parts (i) and (ii))}}

The proof is by induction on $n$, the cases $n=0,1$ being clear by
taking 
$$
v_{()}\;={\bf{0}},\;\;v_{(0)}\;=\;{\bf{0}} + \F_q,\;\;v_{(1)}\;=\;
\vphi {\bf{0}} - \F_q.
$$

Let $n\geq 1$ and consider $\alpha=(\alpha_1,\ldots ,\alpha_t)\in \cE_q(n+1)$.
Set $\beta=(\alpha_1,\ldots ,\alpha_{t-1})$ and $k=N(\beta)$. We have three
cases:

(a) $\alpha_t=0$: We have $v_\beta\in \C[B_q(n)]$. Define
\beq \label{recur1}
v_\alpha = q^kv_\beta + \theta_n(v_\beta) \in \C[B_q(n)]\oplus W(0).
\eeq

(b) $\alpha_t=1$: We have $v_\beta\in \C[B_q(n)]$. Define
\beq \label{recur2}
v_\alpha = \vphi q^{n-k}v_\beta - \theta_n(v_\beta) \in \C[B_q(n)]\oplus W(0).
\eeq

(c) $\alpha_t = \chi,\;\chi\in\cN_q(n)$: We have $v_\beta\in \C[B_q(n-1)]$.
Set $X=X(\chi)$ and define
\beq \label{recur3}
v_\alpha = \lambda(\chi)\mu(X)(v_\beta)\in W(\chi).
\eeq

Let us now check assertions (i) and (ii) in the statement of the theorem,
beginning with (i). We have three cases.

(a) $\alpha_t=0$: Follows from Lemma \ref{kc}. 

(b) $\alpha_t=1$: Follows from Lemma \ref{kc}.

c) $\alpha_t = \chi,\;\chi\in\cN_q(n)$: Set $X=X(\chi)$. It follows from
Theorem \ref{grind1} that
\beqn \cA_{n+1}(v_\alpha)&=&q \;\frac{\vphi q^{n-1-k} - q^k}{q-1}v_\alpha \\
                           &=&\frac{\vphi q^{n+1-(k+1)} - q^{k+1}}{q-1}
                           v_\alpha. 
\eeqn

Assertion (ii) follows from the induction hypothesis using 
the decomposition (\ref{od}), the isomorphism
(\ref{d1}) and observing that the determinant of the $2\times 2$ matrix
\beq \label{2b2det}
&\left[ \ba{cr} q^k & 1 \\ \vphi q^{n-k} & -1 \ea \right]&
\eeq
is nonzero. \eprf

\noi {\bf{Remark}} It is interesting to see where exactly in the proof
above the diagonal entries of $A_n$ are used since 
the main ingredients
of the proof, (\ref{ind}), (\ref{ind1}), {\ref{ind2}), and Theorem
\ref{grind1}, all  hold even when the diagonal entries are 0. 
However, the base case of the induction ($n=0,1$) needs the diagonal
entries. In effect, the presence of the diagonal entries results in the nice
formula for the eigenvalues which makes possible the 
calculations in the proof of Lemma \ref{kc}.

We denote the basis given in part (ii) of 
Theorem \ref{mt} by $\cB_n$
Note that equations (\ref{recur1}), 
(\ref{recur2}), (\ref{recur3}) give an 
inductive procedure to write down $\cB_{n+1}$ 
given $\cB_n$ and $\cB_{n-1}$.
Note also 
that (up to scalars) this basis is canonical in the sense that we have not made
any choices anywhere.

\noi {\bf{Proof of Theorem \ref{mt} (part (iii))}}

\pf The proof is by induction on $n$, the cases $n=0,1$ being clear.

Let $n\geq 1$. We consider two cases:

(i) Let $\beta=(\beta_1,\ldots ,\beta_{t-1})\in\cE_q(n)$.
Set $k=N(\beta)$ and 
\beqn
&\alpha=(\beta_1,\ldots ,\beta_{t-1},0),\;\;\;\; 
\alpha'=(\beta_1,\ldots ,\beta_{t-1},1).
\eeqn 
Given a vectors $u,v\in\C[B_q(n)]$, we shall write $\inp{u}{v}_n$ for
the inner product (\ref{inp}) calculated in $\C[B_q(n)]$ and $\inp{u}{v}_{n+1}$
for the inner product calculated in $\C[B_q(n+1)]$. We have, for
$X\in B_q(n,k)$,
\beqn 
&\inp{X}{X}_n \;=\; \frac{\vphi^k q^{\bin{k}{2}}}{P_q(n)},\;\; 
\inp{X}{X}_{n+1} \;=\; \frac{\vphi^k q^{\bin{k}{2}}}{P_q(n+1)}\;\;
                 \;=\;\frac{1}{1+\vphi q^n}\inp{X}{X}_n,&\\
&\inp{\theta_n(X)}{\theta_n(X)}_{n+1}\;=\;
\frac{\vphi^{k+1} q^{\bin{k+1}{2}}}{P_q(n+1)}\,q^{n-k}
\;=\;\frac{\vphi q^n}{1+\vphi q^n}\inp{X}{X}_n.&
\eeqn
It follows that
\beq \label{inpind}
&\inp{v}{v}_{n+1}\;=\;\frac{1}{1+\vphi q^n}\inp{v}{v}_n,\;\;
\inp{\theta_n(v)}{\theta_n(v)}_{n+1}\;=\;
\frac{\vphi q^n}{1+\vphi q^n}\inp{v}{v}_n,\;\;v\in \C[B_q(n)].&
\eeq
Note that the scalar factors on the right hand side are uniform across 
all vectors
and do not depend on $k$. Thus, since $\C[B_q(n)]$ and $W(0)$ are orthogonal
in $\C[B_q(n+1)]$, 
it follows by the induction hypothesis that
$\{v_\beta , \theta_n(v_\beta)\,|\, \beta\in \cE_q(n)\}$ is an orthogonal basis
of $\C[B_q(n)]\oplus W(0)$.
We have
\beq \label{indfs1}
&v_\alpha\;=\;q^kv_\beta + \theta_n(v_\beta),\;\;
v_{\alpha'}\;=\;\vphi q^{n-k}v_\beta - \theta_n(v_\beta).&
\eeq
Since $v_\beta$ is orthogonal to $\theta_n(v_\beta)$ we have,
using (\ref{inpind}),
\beqn
\inp{v_\alpha}{v_{\alpha'}}_{n+1} &= & 
\vphi q^n\inp{v_\beta}{v_\beta}_{n+1} - \inp{\theta_n(v_\beta)}
{\theta_n(v_\beta)}_{n+1}\\
&=&\frac{\vphi q^n}{1+\vphi q^n}\inp{v_\beta}{v_\beta}_n - 
\frac{\vphi q^n}{1+\vphi q^n}\inp{v_\beta}{v_\beta}_n\\
&=&0.
\eeqn
From the isomorphism $\theta_n$ we now see that
$$ \{ v_\alpha , v_{\alpha'}\;|\; \beta\in \cE_n(q)\}$$
is an orthogonal basis of $\C[B_q(n)]\oplus W(0)$.

(ii) Let $\beta = (\beta_1,\ldots ,\beta_{t-1})\in\cE_q(n-1)$ and 
let $\chi\in\cN_q(n)$.
Set $\alpha=(\beta_1,\ldots ,\beta_{t-1},\chi)\in\cE_q(n+1)$ 
and $X=X(\chi)$, where $X\in B_q(n,n-1)$. We have 
$v_\alpha = \lambda(\chi)\mu(X)(v_\beta)$.

Let $Y\in B_q(X)$ with $\dim(Y)=k$. We have
$$\inp{Y}{Y}_{n-1}\;=\;\frac{\vphi^{k} q^{\bin{k}{2}}}{P_q(n-1)}.
$$ 

Now observe the following: $p(\chi)(\wh{Y})$ is a 
linear combination of the elements of the orbit
$[\wh{Y}]$, whose cardinality is $q^{n-k}$. 
The number of elements $g\in H(n+1,\F_q)$ with $g\cdot \wh{Y}=\wh{Y}$
is $q^k$ and by Lemma \ref{cl} each such $g$ satisfies $\chi(g)=1$.
So, for $Z\in [\wh{Y}]$, if $g_1\cdot \wh{Y} = g_2 \cdot \wh{Y} = Z$ then
$\chi(g_1)=\chi(g_2)$.
Thus we have
\beqn
&\inp{\lambda(\chi)(Y)}{\lambda(\chi)(Y)}_{n+1}\;=\;
q^{-2k}\frac{\vphi^{k+1}q^{\bin{k+1}{2}}}{P_q(n+1)}\,q^{2k}q^{n-k}\;=\;
\frac{\vphi q^n}{(1+\vphi q^{n-1})(1+\vphi q^n)}\inp{Y}{Y}_{n-1}.&
\eeqn
It follows from (\ref{ds}) that
\beq \label{inpind1}
&\inp{\lambda(\chi)\mu(X)(u)}{\lambda(\chi)\mu(X)(v)}_{n+1}\;=\;
\frac{\vphi q^n}{(1+\vphi q^{n-1})(1+\vphi q^n)}\inp{u}{v}_{n-1},
\;\;u,v\in \C[B_q(n-1)].&
\eeq
From the isomorphism $\lambda(\chi)\mu(X)$ we now see that
$$ \{ v_\alpha\;|\; \beta\in \cE_q(n-1)\}$$
is an orthogonal basis of $W(\chi)$.

That completes the proof. \eprf

\section{{Action of $\cA_n^*$ on the eigenspaces of $\cA_n$}}

For $k=0,1,\ldots ,n$ define $E_n(k)\seq \C[B_q(n)]$ to be the eigenspace of
$\cA_n$ with eigenvalue $\frac{\vphi q^{n-k}-q^k}{q-1}$. In this section we
consider the action of $\cA_n^*$ on the eigenspaces of $\cA_n$ and give a
new proof of Theorem \ref{qp}.

Any homogeneous subspace of $\C[B_q(n)]$ is closed under the action of 
$\cA_n^*$ and so the subspaces $\C[B_q(n)]$, $W(0)$, $W(\chi),\;\chi\in
\cN_q(n)$ of $\C[B_q(n+1)]$ are all closed under $\cA_{n+1}^*$.

Note the following relations in $\C[B_q(n)]\oplus W(0)$:
\beqn  
\cA_{n+1}^*(X) &=& \cA_n^*(X),\;\;X\in B_q(n),\\
\cA_{n+1}^*(\theta_n(X)) &=& q^{-1}\theta_n (\cA_n^*(X)),\;\;X\in B_q(n).\\
\eeqn
It follows that
\beq \label{dind1}  
\cA_{n+1}^*(v) &=& \cA_n^*(v),\;\;v\in \C[B_q(n)],\\ 
\label{dind2}
\cA_{n+1}^*(\theta_n(v)) &=& q^{-1}\theta_n (\cA_n^*(v)),\;\;v\in \C[B_q(n)].
\eeq

Let $\chi \in \cN_q(n)$ and $X=X(\chi)$. Then we have
\beqn
\cA_{n+1}^*(\lambda(\chi)\mu(X)(Y))=q^{-1}
\lambda(\chi)\mu(X)(\cA_{n-1}^*(Y)),\;\;Y\in B_q(n-1).
\eeqn
It follows that
\beq \label{dind3}
\cA_{n+1}^*(\lambda(\chi)\mu(X)(v))=q^{-1}
\lambda(\chi)\mu(X)(\cA_{n-1}^*(v)),\;\; v\in \C[B_q(n-1)].
\eeq
We can now give the 

\noi 
{\bf{Proof of Theorem \ref{qp}}} 
The proof is by induction on $n$, the cases
$n=0,1$ being clear.

Consider a basis element $v_\beta \in \cB_{n+1}$ with $N(\beta)=k$. Then
$v_\beta\in E_{n+1}(k)$. There are three cases to consider:

(i) The last element of the sequence $\beta$ is 0: 
We have $v_\beta = q^k v_\alpha + \theta_n(v_\alpha)$, where
$v_\alpha \in \cB_{n}$ with $N(\alpha)=k$. Then $v_\alpha \in E_n(k)$.

By the induction hypothesis we have 
\beqn
\cA_n^*(v_\alpha)&=& \sum_{j=k-1}^{k+1} u_j,\;\;u_j\in E_n(j),\;j=k-1,k,k+1.
\eeqn
We have, by (\ref{dind1}) and (\ref{dind2}),
\beqn
\cA_{n+1}^*(v_\beta)&=& \cA_{n+1}^*(q^k v_\alpha + \theta_n(v_\alpha))\\
                    &=& q^k \cA_n^*(v_\alpha) + q^{-1}
\theta_n(\cA_n^*(v_\alpha))\\
&=& q^{k} \left\{ \sum_{j=k-1}^{k+1} u_j\right\} + q^{-1} \left\{
\sum_{j=k-1}^{k+1} \theta_n(u_j)\right\}. 
\eeqn
We consider the three terms on the right hand side separately.

(a) We have
\beqn
q^k u_{k+1} + q^{-1}\theta_n(u_{k+1})&=&q^{-1}  ( q^{k+1} u_{k+1} +
\theta_n(u_{k+1})) \in E_{n+1}(k+1),
\eeqn
by Lemma \ref{kc} (since $u_{k+1}\in E_n(k+1)$).

(b) We have
\beqn
q^k u_{k} + q^{-1}\theta_n(u_{k})&=&q^{-1}  ( q^{k+1} u_{k} +
\theta_n(u_{k})). 
\eeqn
By the nonsingularity of the $2\times 2$ matrix (\ref{2b2det}), 
the vector above is in the span of the vectors
$q^k u_k + \theta_n(u_k) \in E_{n+1}(k)$ and
$\vphi q^{n-k}u_k - \theta_n(u_k) \in E_{n+1}(k+1)$.

(c) We have
\beqn
q^{k} u_{k-1} + q^{-1}\theta_n(u_{k-1})&=&q^{-1}  ( q^{k+1} u_{k-1} +
\theta_n(u_{k-1})). 
\eeqn
Similar to the subcase (b) above, the vector above is in the span of the vectors
$q^{k-1} u_{k-1} + \theta_n(u_{k-1}) \in E_{n+1}(k-1)$ and
$\vphi q^{n-k+1}u_{k-1} - \theta_n(u_{k-1}) \in E_{n+1}(k)$.

That completes the proof of the first case.

(ii) The last element of the sequence $\beta$ is 1:
We have $v_\beta = \vphi q^{n-k+1} v_\alpha - \theta_n(v_\alpha)$, where
$v_\alpha \in \cB_{n}$ with $N(\alpha)=k-1$. Then $v_\alpha \in E_n(k-1)$.

By the induction hypothesis we have 
\beqn
\cA_n^*(v_\alpha)&=& \sum_{j=k-2}^{k} w_j,\;\;w_j\in E_n(j),\;j=k-2,k-1,k.
\eeqn
We have, by (\ref{dind1}) and (\ref{dind2}),
\beqn
\cA_{n+1}^*(v_\beta)&=& \cA_{n+1}^*(\vphi q^{n-k+1} 
v_\alpha - \theta_n(v_\alpha))\\
                    &=& \vphi q^{n-k+1} \cA_n^*(v_\alpha) - q^{-1}
\theta_n(\cA_n^*(v_\alpha))\\
&=& \vphi q^{n-k+1} \left\{ \sum_{j=k-2}^{k} w_j\right\} - q^{-1} \left\{
\sum_{j=k-2}^{k} \theta_n(w_j)\right\}. 
\eeqn

We consider the three terms on the right hand side separately.

(a) We have
\beqn
\vphi q^{n-k+1} w_{k-2} - q^{-1}\theta_n(w_{k-2})&=&
q^{-1}  ( \vphi q^{n-k+2} w_{k-2} - 
\theta_n(w_{k-2})) \in E_{n+1}(k-1),
\eeqn
by Lemma \ref{kc} (since $w_{k-2}\in E_n(k-2)$).

(b) We have
\beqn
\vphi q^{n-k+1} w_{k-1} - q^{-1}\theta_n(w_{k-1})&=&
q^{-1}  ( \vphi q^{n-k+2} w_{k-1} - 
\theta_n(w_{k-1})).
\eeqn
The vector above is in the span of the vectors
$q^{k-1} w_{k-1} + \theta_n(w_{k-1}) \in E_{n+1}(k-1)$ and
$\vphi q^{n-k+1}w_{k-1} - \theta_n(w_{k-1}) \in E_{n+1}(k)$.

(c) We have
\beqn
\vphi q^{n-k+1} w_{k} - q^{-1}\theta_n(w_{k})&=&
q^{-1}  ( \vphi q^{n-k+2} w_{k} - 
\theta_n(w_{k})).
\eeqn
The vector above is in the span of the vectors
$q^{k} w_{k} + \theta_n(w_{k}) \in E_{n+1}(k)$ and
$\vphi q^{n-k}w_k - \theta_n(w_{k}) \in E_{n+1}(k+1)$.

That completes the proof of the second case.

(iii) The last element of the sequence $\beta$ is $\chi$, for some
$\chi\in \cN_q(n)$: Let $X=X(\chi)$. Then
$v_\beta = \lambda(\chi)\mu(X)(v_\alpha)\in W(\chi)$, where $N(\alpha)=k-1$
and $v_\alpha \in E_{n-1}(k-1)$. The result now follows from Theorem
\ref{grind1} and (\ref{dind3}). 
\eprf

\begin{center}
\section{{Acknowledgement}}
\end{center}

This work was motivated by e-mail discussions with 
Paul Terwilliger. I am very grateful to him
for his encouragement and support. I thank Subhajit
Ghosh for useful feedback and pointing out {\bf\cite{kk}}.

\end{document}